\documentclass[letterpaper, 10 pt, conference]{ieeeconf}  

\IEEEoverridecommandlockouts                              
\usepackage{graphicx} 
\usepackage{amsmath} 
\usepackage{amssymb}  
\usepackage{cite}

\title{\LARGE \bf
Interpreting Reinforcement Learning Model Behavior via Koopman with Control
}

\author{William T. Redman$^1$ 
\thanks{$^1$W.T. Redman is with the Electrical and Computer Engineering Department, Johns Hopkins University and AIMdyn Inc. Correspondence: {\tt\small wredman4@jh.edu}}
}

\begin{document}

\maketitle
\thispagestyle{empty}
\pagestyle{empty}

\begin{abstract}
Reinforcement learning (RL) models have shown the capability of learning complex behaviors, but quantitatively assessing those behaviors -- which is critical for safety assurance and the discovery of novel strategies -- is challenging. By viewing RL models as control systems, we hypothesize that data-driven approximations of their associated Koopman operators may provide dynamical information about their behavior, thus enabling greater interpretability. To test this, we apply the Koopman with control framework to RL models trained on several standard benchmark environments and demonstrate that properties of the fit linear control models, such as stability and controllability, evolve during training in a task dependent manner. Comparing these metrics across different training epochs or across differently optimized RL models enables an understanding of how they differ. In addition, we find cases where -- even when the reward achieved by the RL model is static -- the stability and controllability is nonetheless evolving, predicting increased reward with further training. This suggests that these metrics may be able to serve as hidden progress measures, a core idea in mechanistic interpretability. Taken together, our results illustrate that the Koopman with control framework provides a comprehensive way in which to analyze and interpret the behavior of RL models, particularly across training. 

\end{abstract}

\section{INTRODUCTION}

Reinforcement learning (RL) models have exhibited remarkable success on a wide range of applied problems \cite{mnih2015human, silver2018general, degrave2022magnetic}. However, in many cases, how they achieve this high performance remains unknown. Attempts to quantitatively interpret the behaviors learned by RL models can require expert knowledge \cite{mcgrath2022acquisition}, limiting the ability to achieve safety assurance and to identify novel capabilities that are discovered by RL models during training \cite{schut2025bridging}. Even in simple environments, understanding how the behavior of RL models changes with learning is non-trivial, with typical approaches using heuristics to identify changes in the distribution of visited states.

\begin{figure}
    \centering
    \includegraphics[width=0.90\linewidth]{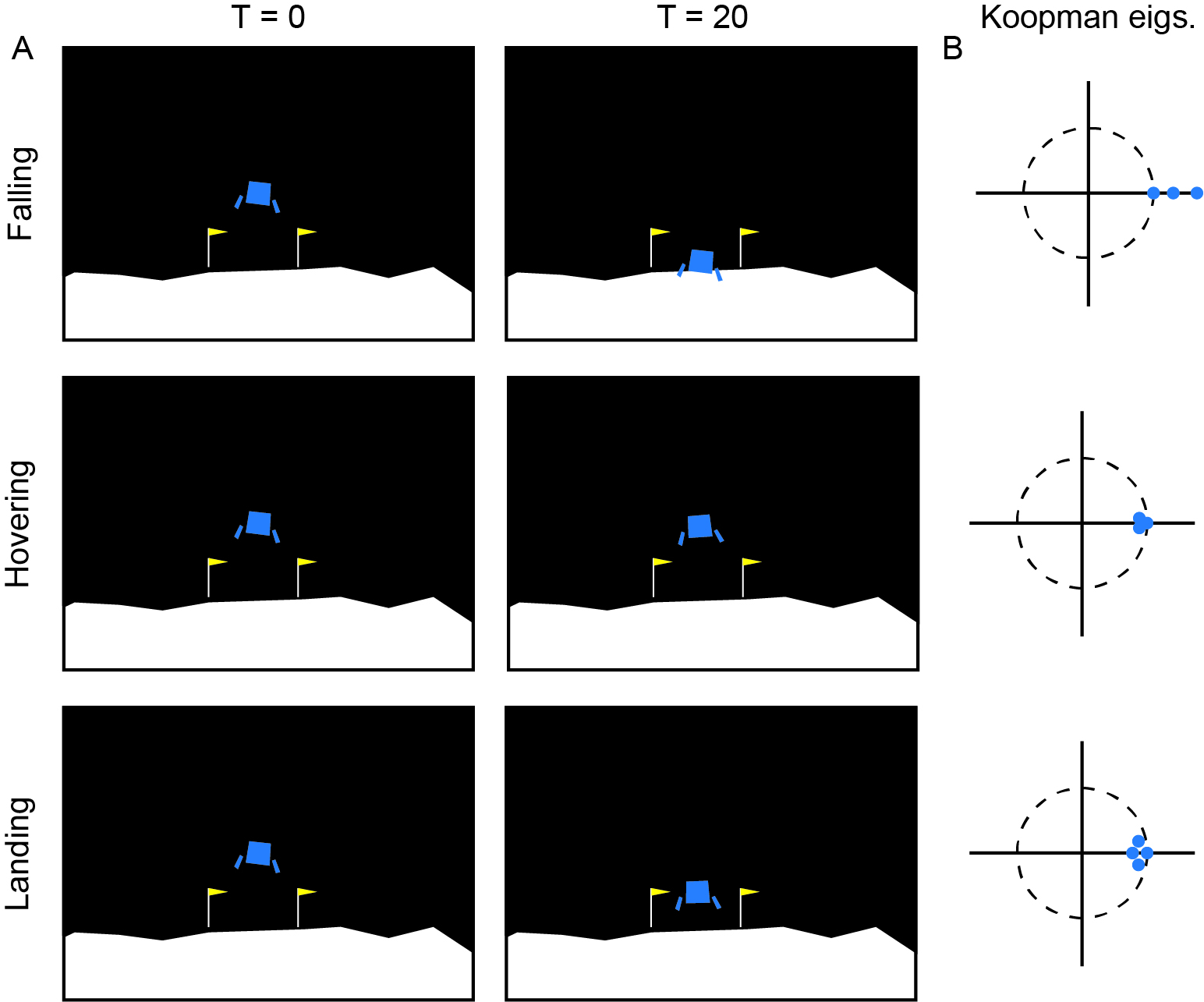}
    \caption{\textbf{Schematic illustration of how different behaviors of RL models may have different associated Koopman eigenvalues.} (A) Three example dynamical behaviors (``falling'', ``hovering'', ``landing'') that are possible in the \texttt{LunarLander} RL environment. First column is at the start of the trial ($T = 0$) and second column is $20$ time steps later ($T = 20$). (B) Possible Koopman eigenvalues associated with each of the behaviors in (A).}
    \label{fig: Koopm4RL Schematic}
\end{figure}

At their core, RL models are control systems, taking in observations about the environment and performing actions (inputs), which lead to new states (outputs). While the state transition probabilities may be fixed in a given environment (e.g., fixed by the laws of physics, the rules of chess), the range of dynamical behaviors performed by an RL model (Fig. \ref{fig: Koopm4RL Schematic}A: ``falling'', ``hovering'', ``landing'') may change across training and/or between differently optimized RL models. A natural way to provide better interpretability of RL model behavior would therefore be to leverage ideas from control theory, such as controllability and stability \cite{hespanha2018linear}. However, because of the complexity and nonlinearity of the environments RL models are applied in, directly utilizing tools from control theory to provide insight into their behavior can be challenging. 

Koopman operator theory \cite{mezic2005spectral, budivsic2012applied}, a framework for learning linear representations of nonlinear dynamical systems by lifting the underlying state-space to a high-dimensional function space, has shown considerable potential for providing interpretable models of complex, real-world systems \cite{rowley2009spectral, avila2020data}, including algorithms \cite{dietrich2020koopman, redman2022algorithmic} and machine learning models \cite{ostrow2023beyond, akrout2023representations, redman2024identifying, huang2025inputdsa}. Extensions of Koopman operator theory to include control have been developed \cite{korda2018linear, proctor2018generalizing}, enabling the use of linear control for highly nonlinear systems \cite{bruder2020data, haggerty2023control}. Recent work has shown that such Koopman operator theory with control can be used to compare the internal activation dynamics of recurrent neural networks (RNNs) trained with RL \cite{huang2025inputdsa}. However, the extent to which Koopman with control can broadly enhance the interpretability of RL model behavior has yet to be explored. 

Here, we show that linear time-invariant (LTI) control models, of the form $z_{t + 1} = Az_t + Bu_t$ (where $z$ is a lifted version of the state-space variables $x$), can be fit to the output trajectories generated by individual RL models, and that properties of $A$ and $B$ can provide insight into the associated RL model's behavior. We demonstrate that, across training, these metrics evolve in a task dependent manner and are often, although not always, aligned with the RL model's reward. This enables a quantitative understanding of the way in which the RL models behavior evolves and can be used to compare RL models trained with different optimizers (e.g., PPO \cite{schulman2017proximal}, A2C \cite{mnih2016asynchronous}). In some settings, we find that these metrics can exhibit changes across training, even when the reward appears static. This suggests that the properties of $A$ and $B$ can act as hidden progress measures \cite{barak2022hidden}, an important concept in the study of the mechanistic interpretability \cite{nanda2023progress}.

\section{KOOPMAN OPERATOR THEORY WITH CONTROL FOR STUDYING RL BEHAVIOR}
\label{section: Koopman operator theory with control}

\subsection{RL models as control systems}
\label{subsection: RL as control systems}
RL models take in observations, $o \in \mathcal{O}$, and perform actions, $a \in \mathcal{A}$, which evolve the underlying state of the system, $s \in \mathcal{S}$. Thus, RL models are defined on tuples $(\mathcal{O}, \mathcal{S}, \mathcal{A})$, where $\mathcal{O}$ is the observation space, $\mathcal{S}$ is the state-space, and $\mathcal{A}$ is the action space. By viewing the actions as inputs and the observations as outputs, we can recognize the RL model update 
\begin{equation}
    \label{eq: RL model update}
    \begin{split}
        s_{t + 1} = g(s_t, a_t) \\
        o_{t} = h(s_t, a_t) \\ 
    \end{split}
\end{equation}
as a nonlinear control system, where $a_t = m(o_t)$ is the action generated by the RL model, through learning the policy $m: \mathcal{O} \rightarrow \mathcal{A}$. Here, $g$ and $h$ are maps defined by the environment, where $g$ defines how a given action translates to a transition in states (e.g., how applying a specific torque to a pendulum affects its angular position) and $h$ defines what the RL model is able to observe of state $s_t$. In the case where there is full observability, $o_t = s_t$. In general, modern RL models make actions probabilistically, meaning the policy $m$ is not deterministic. In some complex environments, actions can lead to probabilistic transitions between states, meaning $g$ is also not deterministic. Furthermore, $g$ may be time-varying, making the RL problem (or correspondingly, the control problem) considerably more challenging. In this work, we consider environments where $g$ is deterministic and time-invariant, and where the states are fully observable (i.e., $h = \text{Id}_{s_t}$). Whether and how our approach can be extended to more difficult environments is an important future direction. 

\subsection{Koopman operator theory}
\label{subsection: Koopman operator theory}

While a discrete-time dynamical system, $x_{t + 1} = T(x_t)$ with $T: \mathcal{X} \rightarrow \mathcal{X}$, may be governed by a nonlinear map $T$, a linear representation may be achieved by lifting the states of the system, $x \in \mathcal{X}$, to an infinite dimensional function space, $\mathcal{F}$. Functions in this space (referred to as observables), $f \in \mathcal{F}$, evolve via the action of a composition operator, $U$,
\begin{equation}
\label{eq: Koopman operator}
    U^t f(x_0) = f[T^t(x_0)].
\end{equation}
The operator $U$ is referred to as the Koopman operator and finite approximations of it, via data-driven numerical methods  \cite{rowley2009spectral, williams2015data, arbabi2017ergodic, mauroy2024analytic, colbrook2024rigorous}, provide linear models of complex dynamical systems. The linearity enables a mode decomposition \cite{mezic2005spectral}\footnote{In general, there is an additional term in Eq. \ref{eq: Koopman mode decomposition} corresponding to the continuous part of the spectrum. Given that, in order to be performant, RL models are expected not to be chaotic, we assume that the systems we fit only have point spectra.}, 
\begin{equation}
    \label{eq: Koopman mode decomposition}
    U^t f(x_0) = \sum_{k = 1}^N \lambda_k^t \phi_k(x_0) v_k,
\end{equation}
where $(\lambda_k, \phi_k, v_k)$ are the $N$ Koopman eigenvalues and their corresponding eigenfunctions and Koopman modes, respectively. This decomposition enables greater interpretability of the complex dynamics by extracting temporal and spatial information, which can be used to gain insight into the spatio-temporal properties of the system \cite{rowley2009spectral, avila2020data}.

\subsection{Koopman with control}
\label{subsection: Koopman with control}

Extensions of Koopman operator theory to include systems with control have been developed \cite{korda2018linear, proctor2018generalizing}. In particular, for a control system described by
\begin{equation}
    \label{eq: control systems}
    \begin{split}
    x_{t + 1} = g(x_t, u_t) \\
    y_t = h(x_t, u_t)
    \end{split}
\end{equation}
where $x_t$ is the state, $u_t$ is the input, and $y_t$ is the output at time $t$, an LTI surrogate model can be learned, 
\begin{equation}
    \label{eq: Koopman with control}
    \begin{split}
        z_{t + 1} = A z_t + B u_t \\
        y_t = C z_t + D u_t,
    \end{split}
\end{equation}
where $z_t$ is the lifted version of the state $x_t$ (i.e., $z_t = [f_1(x_t), ..., f_m(x_t)]$, for some $\{f_i\} \in \mathcal{F}$). 

\subsection{Koopman with control for studying RL behavior}
\label{subsection: Koopman with control for RL}

While the map $g: \mathcal{S} \rightarrow \mathcal{S}$ is defined by the environment (e.g., laws of physics, rules of chess), different dynamical behaviors can demonstrate different aspects of $g$. For instance, in the standard RL environment \texttt{LunarLander}, an RL model controls a simplified lunar lander and is tasked with landing it on the surface of the moon, within a specified spatial window (Fig. \ref{fig: Koopm4RL Schematic}A -- yellow flags denote target location of landing). To succeed, it must stabilize itself by coordinating its left and right engines. In addition, it must slow itself on its descent, so it does not crash (Fig. \ref{fig: Koopm4RL Schematic}A, bottom row). Fitting a Koopman with control model (Eq. \ref{eq: Koopman with control}) to trajectories where the RL model has successfully learned this behavior, the Koopman eigenvalues associated with the matrix $A$ will have norm close to, but less than, $1$ (Fig. \ref{fig: Koopm4RL Schematic}B, bottom row). This corresponds to the slow and stable descent of the well trained lunar lander. In contrast, if the RL model has not been fully trained, it may exhibit trajectories where the lunar lander falls (Fig. \ref{fig: Koopm4RL Schematic}A, top row) or where the lunar lander only hovers (Fig. \ref{fig: Koopm4RL Schematic}A, middle row). Each of these behaviors will have different associated Koopman eigenvalues (Fig. \ref{fig: Koopm4RL Schematic}B). This suggests that properties of $A$, and potentially $B$, can capture information that may prove insightful for understanding RL model behavior.  

Here, we consider two such properties of the fit $A$ and $B$ matrices: maximum eigenvalue norm of $A$ and rank of the controllability matrix\footnote{To make the rank more meaningful, we normalize so that a value of $1$ denotes that the controllability matrix is full rank.} ($[B, AB, A^2B, ..., A^{n - 1}B]$). These capture aspects of the RL model's stability and controllability. Other properties of the learned Koopman with control model may provide further insight and should be considered in future work. 

\section{RESULTS}

To demonstrate the potential of utilizing Koopman with control to interpret the behavior of RL models, we consider three standard RL environments: \texttt{CartPole}, \texttt{Acrobot}, and \texttt{LunarLander}. We train RL models using PPO \cite{schulman2017proximal} and A2C \cite{mnih2016asynchronous}. We implement our experiments using Gymnasium \cite{towers2024gymnasium} and Stable-Baselines3 \cite{raffin2021stable}. To fit the Koopman with control models, we use the DMDc algorithm \cite{proctor2018generalizing}, implemented by the PyDMD package \cite{demo2018pydmd, ichinaga2024pydmd}. In all cases, we ensure that the fit model has low mean-squared reconstruction error ($< 0.01$). Code is publicly available\footnote{ $\texttt{https://github.com/Dynamical-Intelligence-Group/} \\ \texttt{Koopman-with-control-for-RL-model-behavior}$}. 

\subsection{Cart Pole} 
\label{subsection: Cart Pole}

\begin{figure}
    \centering
    \includegraphics[width=0.95\linewidth]{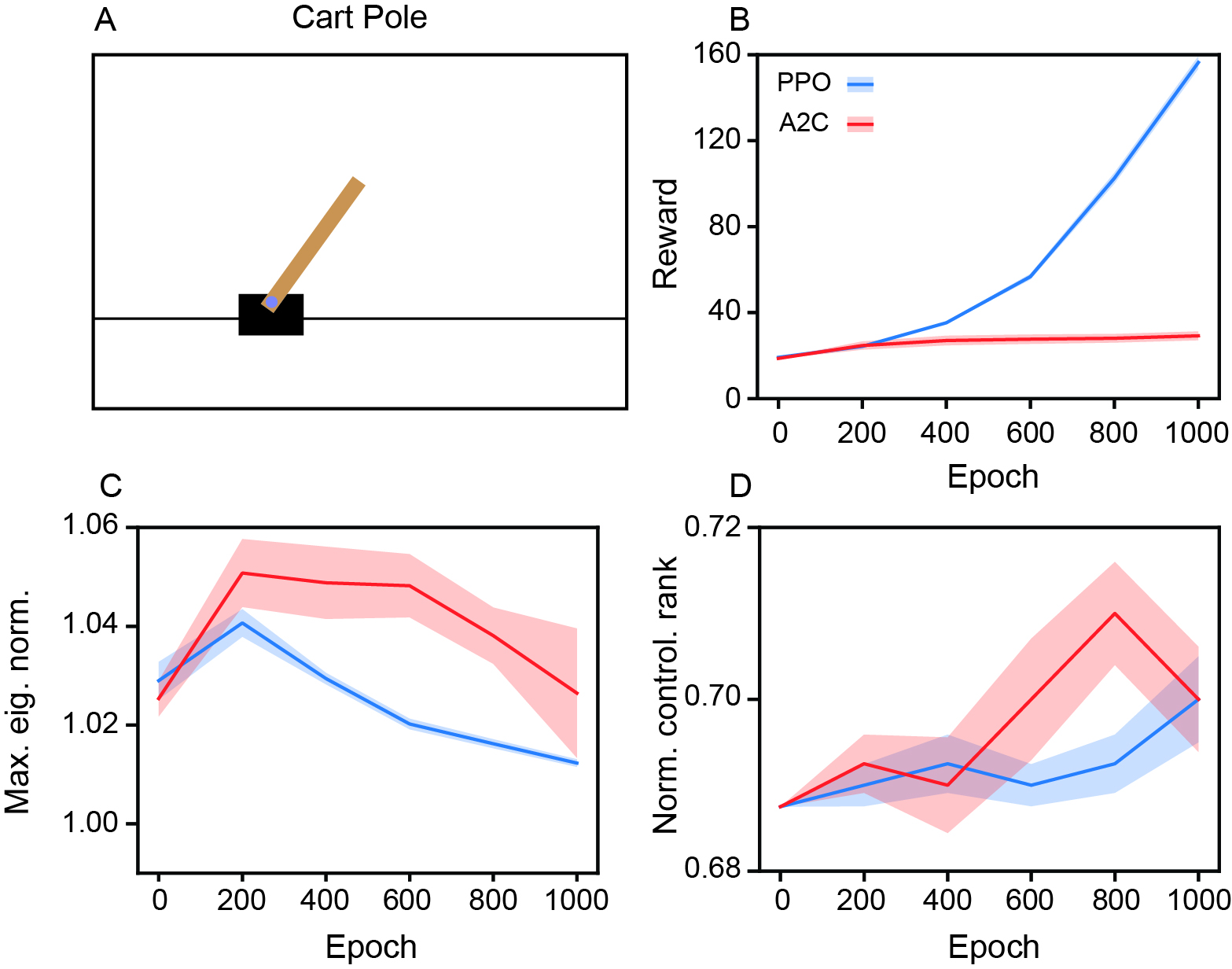}
    \caption{\textbf{Performant RL models exhibit more stable trajectories, with greater controllability, across training on \texttt{CartPole}.} (A) Schematic illustration of the \texttt{CartPole} environment. (B) Median reward, as a function of training epoch, for RL models optimized using PPO (blue) and A2C (red). (C) Maximum eigenvalue norm, as a function of epoch. (D) Normalized rank of controllability matrix, as a function of epoch. (B)--(D) Solid line represents mean and shaded area represents $\pm$ standard error across 25 independently trained RL models.} 
    \label{fig: Cart Pole}
\end{figure}

\texttt{CartPole} is a classic control task \cite{barto1983neuronlike}, in which a pole is attached to a cart via an un-actuated joint (Fig. \ref{fig: Cart Pole}A) and the RL model is trained to stabilize the pole by moving either to the left or the right. These two actions are represented as $0$ and $1$, respectively. This makes the action space $\mathcal{A} = \{0, 1\}$. Because having the input $0$ denote a movement that is equal and opposite to $1$ may make the learning of the associated Koopman with control model more challenging, we consider a ``one-hot'' embedding of the actions, with $u = [1, 0]^{^{\intercal}}$ denoting a movement of the cart to the left and $u = [0, 1]^{^{\intercal}}$ denoting a movement of the cart to the right. The state-space is four-dimensional, comprising of the linear position of the cart, the linear velocity of the cart, the angular position of the pole, and the angular velocity of the pole. For every time-step that the pole is sufficiently upright and the cart is within a fixed spatial interval, the reward is increased by $+1$.  

We train RL models for $1000$ epochs, using either PPO or A2C. Every $200$ epochs, we save the RL model and sample $100$ new trials, each at most $200$ time-steps long. We record the states and actions across all these trials and construct a Koopman with control model that fits the dynamics using time-delay embeddings of the state \cite{brunton2017chaos, arbabi2017ergodic}. We use $n_\text{delay} = 4$ time-delays and fix the SVD rank of the DMDc computation to be $0.95$. We perform this process on $25$ independently initialized and trained RL models to check the robustness of our results. 

We find that, after the first $200$ training epochs, the PPO optimized RL models outperform those optimized with A2C (Fig. \ref{fig: Cart Pole}B). Examining the properties of the fit Koopman with control models, we find that the PPO trained RL models exhibit trajectories with greater stability (maximum eigenvalue norm closer to $1$) (Fig. \ref{fig: Cart Pole}C). The normalized ranks of the controllability matrix (Fig. \ref{fig: Cart Pole}D) are similar, in both cases increasing with training.

While the reward of the A2C optimized RL models is relatively static during the first $1000$ epochs (Fig. \ref{fig: Cart Pole}B), closer examination at the later part of training shows an increase in stability (decreasing maximum eigenvalue norm -- Fig. \ref{fig: Cart Pole}C) and an increase in controllability (increasing normalized rank of controllability matrix -- Fig. \ref{fig: Cart Pole}D). This suggests that, while these models may not be showing much change -- when looking just at the reward -- they may nonetheless be improving.  

\begin{figure}
    \centering
    \includegraphics[width=0.999\linewidth]{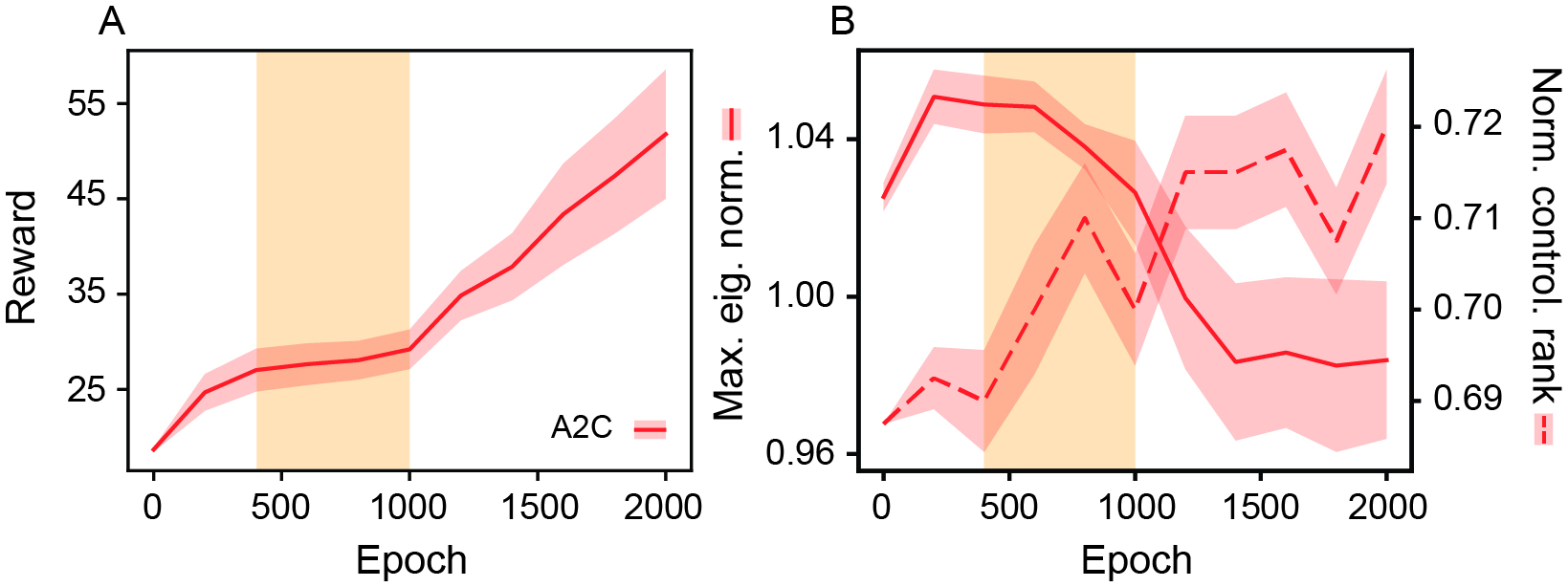}
    \caption{\textbf{Koopman eigenvalue norm and controllability matrix rank contain information that can act as a hidden progress measure.} (A) Median reward, as a function of training epoch, for RL models optimized using A2C. Orange shaded area denotes training time window where little change is seen in the reward. (B) Maximum eigenvalue norm and normalized controllability matrix rank, as a function of epoch. Orange shaded area denotes same training time window as in (A). Note that, for these two metrics, there is an increase in stability and controllability, during this training time window. }
    \label{fig: Cart Pole Hidden Progress Measure}
\end{figure}

To investigate the hypothesis, we train RL models using A2C for another $1000$ epochs. Consistent with increasing stability and controllability (Fig. \ref{fig: Cart Pole Hidden Progress Measure}B, orange shaded area), we find that the reward begins to increase after $1000$ epochs and sees a large improvement in performance by $1600$ epochs (Fig. \ref{fig: Cart Pole Hidden Progress Measure}A). This suggests that the properties of the fit $A$ and $B$ matrices in the Koopman with control model may serve as hidden progress measures \cite{barak2022hidden}, identifying changes in the behavior that are not captured by the reward. 

\subsection{Acrobot}
\label{subsection: Acrobot}

\begin{figure}
    \centering
    \includegraphics[width=0.95\linewidth]{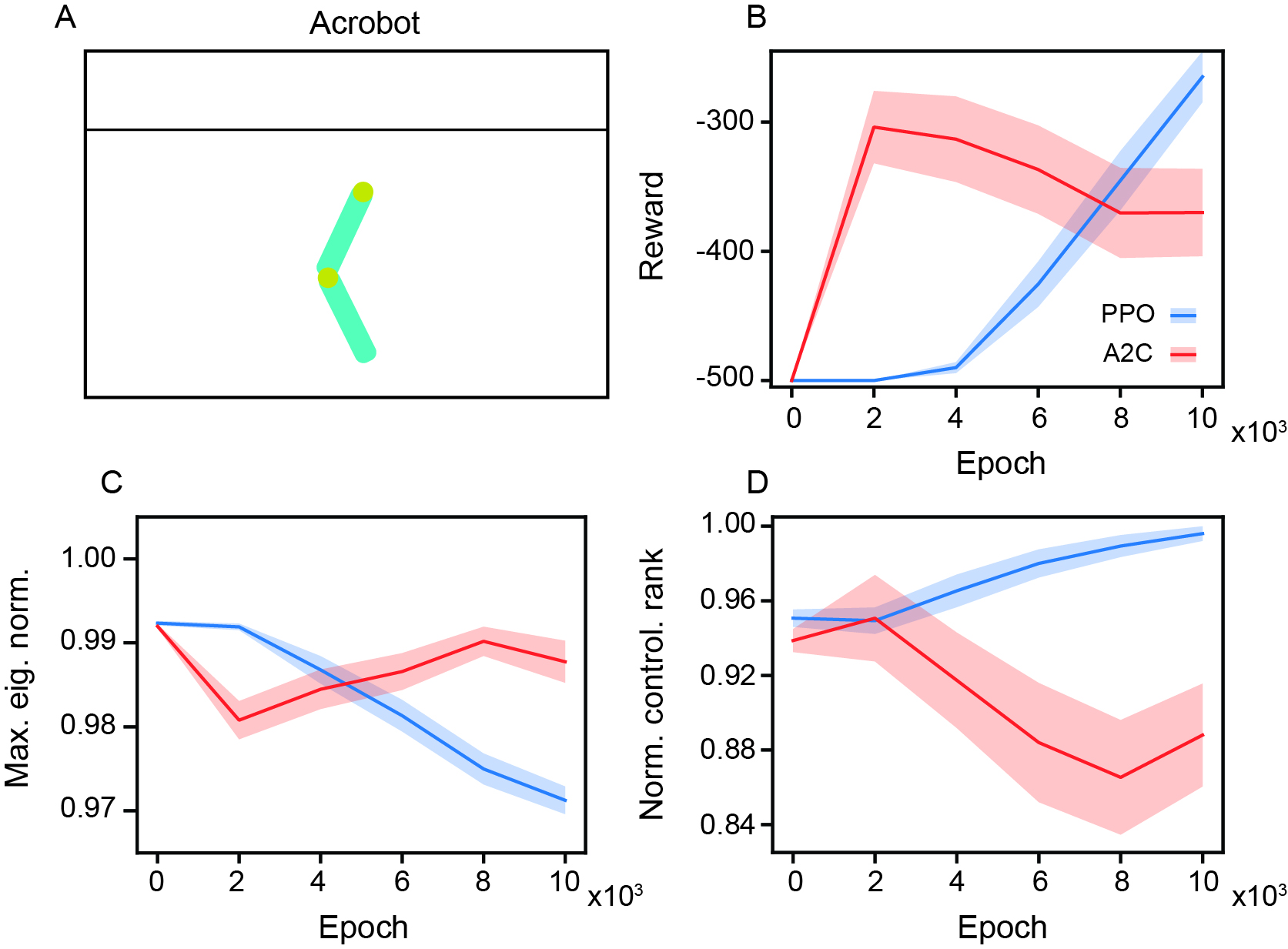}
    \caption{\textbf{Performant RL models exhibit less stable trajectories, with a greater controllability, across training on \texttt{Acrobot}.} (A) Schematic illustration of the \texttt{Acrobot} environment. (B) Reward, as a function of training epoch, for RL models optimized using PPO (blue) and A2C (red). (C) Maximum norm of eigenvalues, as a function of training epoch. (D) Normalized rank of controllability matrix, as a function of epoch. (B)--(D) Solid line represents mean and shaded area represents $\pm$ standard error across 25 independently trained RL models.} 
    \label{fig: Acrobot}
\end{figure}

Another standard classic control RL environment is \texttt{Acrobot} \cite{sutton1995generalization}. In this environment, an RL model must control two links that are linearly connected with an actuated joint (Fig. \ref{fig: Acrobot}A). One end of the chain is fixed and RL models are trained to provide torque on the actuated joint in a sufficient way so as to get the free end of the chain above a given height (Fig. \ref{fig: Acrobot}A, black line). The action-space is given by $\mathcal{A} = \{-1, 0, 1\}$, where $a \in \mathcal{A}$ corresponds to applying $a$ torque to the actuated joint. Because of the success in Sec. \ref{subsection: Cart Pole}, we again use a one-hot embedding of the action to model the control. The state-space comprises of the cosine and sine of $\theta_1$ and $\theta_2$, where $\theta_1$ is the angle of the first link and $\theta_2$ is the relative angle between the two links. In addition, the state-space includes the angular velocity of $\theta_1$ and $\theta_2$, making $|\mathcal{S}| = 6$. For every time-step that the free end of the chain is below the goal height, the reward is decreased by $-1$. 

We train RL models for $10^4$ epochs, using PPO and A2C. Every $2000$ epochs, we save the RL model and sample $100$ new trials, each at most $500$ time-steps long. We use $n_\text{delay} = 5$ time-delays and fix the SVD rank of the DMDc computation to be $0.99$. We perform this process on $25$ independently initialized and trained RL models.

We find that the A2C trained models quickly outperform the PPO trained models, reaching a reward of approximately $-300$ within $2000$ epochs (Fig. \ref{fig: Acrobot}B). The better reward achieved by RL models trained using A2C is accompanied by a rapid decrease in maximum eigenvalue norm within the first $2000$ training epochs (Fig. \ref{fig: Acrobot}C), reflecting the need for transient (as opposed to stable) behavior to launch the free end of the chain above the target height. While the RL models trained with PPO are less performant initially, there is a steady improvement in reward after the first $2000$ training epochs, with these models ultimately outperforming the models trained with A2C (Fig. \ref{fig: Acrobot}B). This is accompanied by a larger decrease in maximum eigenvalue norm (Fig. \ref{fig: Acrobot}C) and  larger increase in normalized controllability rank (Fig. \ref{fig: Acrobot}D.

\subsection{Lunar Lander}
\label{subsection: Lunar Lander}

\begin{figure}
    \centering
    \includegraphics[width=0.95\linewidth]{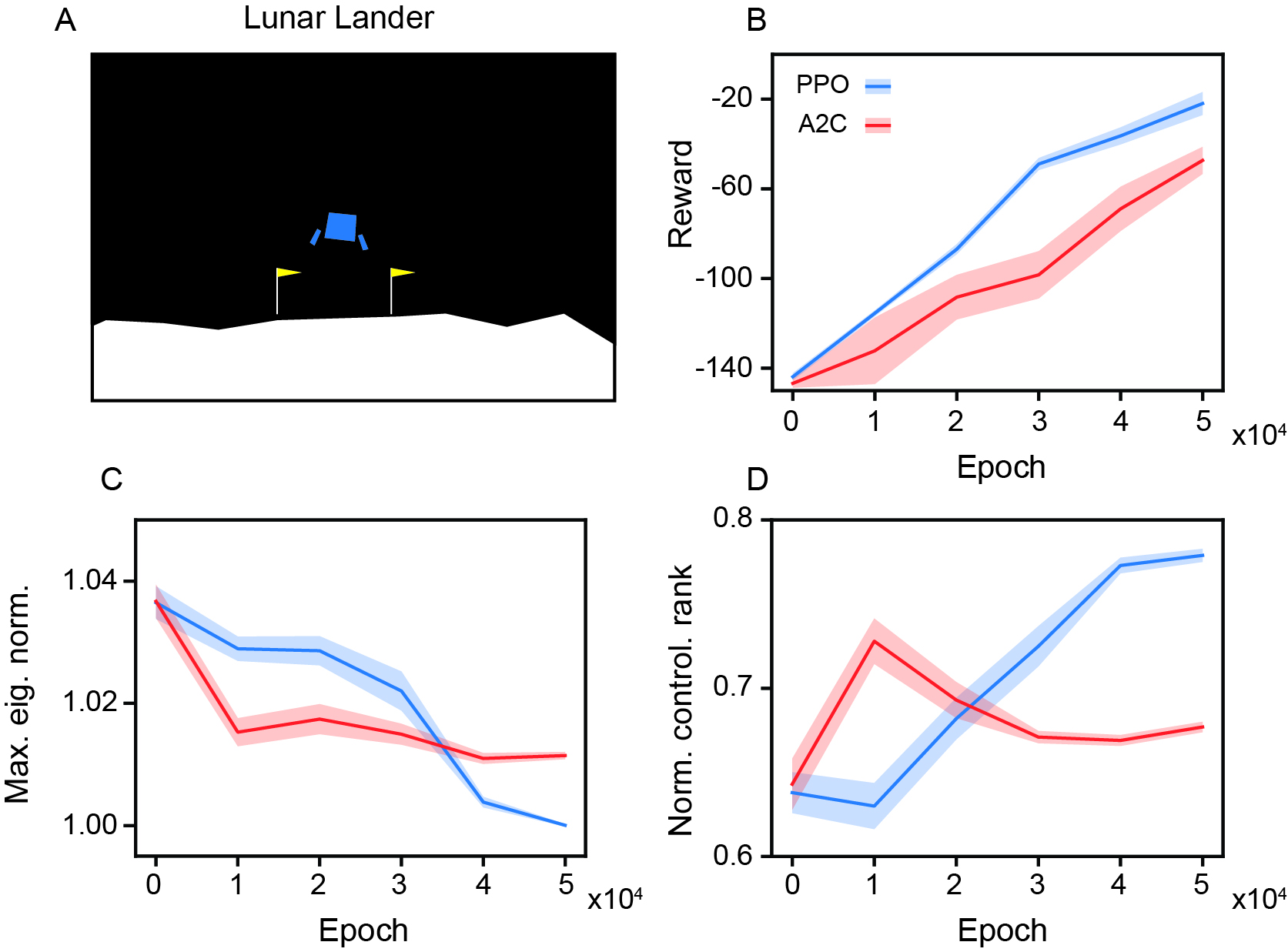}
    \caption{\textbf{Performant RL models exhibit more stable trajectories, with greater controllability, across training on \texttt{LunarLander}.} (A) Schematic illustration of the \texttt{LunarLander} environment. (B) Median reward, as a function of training epoch, for RL models optimized using PPO (blue) and A2C (red). (C) Maximum norm of eigenvalues, as a function of epoch. (D) Normalized ank of controllability matrix, as a function of epoch. (B)--(D) Solid line represents mean and shaded area represents $\pm$ standard error across 25 independently trained RL models. } 
    \label{fig: Lunar Lander}
\end{figure}

Finally, we apply the Koopman with control framework to study the behavior of RL models trained in \texttt{LunarLander}, a more complex control problem that is widely used as a benchmark in RL. In this environment, an RL model is trained to guide the lunar lander to a specified spatial range of the moon (Fig. \ref{fig: Lunar Lander}A). The RL model is able to make one of four actions at each time-step: do nothing, fire the left orientation engine, fire the main engine, fire the right orientation engine. These are denoted as $0$ to $3$, respectively, making $\mathcal{A} = \{0, 1, 2, 3\}$. As in Sec. \ref{subsection: Cart Pole}, we believe this representation of actions may be not ideal for the Koopman with control model to properly represent the way in which the actions affect the lifted states. Therefore, we again use a one-hot embedding of the actions to model the control inputs. The state-space consists of the $x-$ and $y-$coordinates of the lunar lander, the linear velocities along the $x$ and $y$ axes, the angle of the lunar lander, the angular velocity, and two scalars that are $0$ if the left (right) leg of the lander has not contacted the moon's surface and $1$ otherwise. Thus, $|\mathcal{S}| = 8$. Unlike \texttt{CartPole} and \texttt{Acrobot}, the reward function contains multiple terms, including a bonus for landing ($+100$), a penalty for crashing ($-100$), increases (decreases) for slower (faster) descent, and a penalty for firing the engines. 

We train RL models for $5 \cdot 10^4$ epochs, using PPO and A2C as optimizers. Every $10^4$ epochs, we save the RL model and sample $100$ new trials, each at most $1000$ time-steps long. We use $n_\text{delay} = 5$ time-delays and fix the SVD rank of the DMDc computation to be $0.99$. We perform this process on $25$ independently initialized and trained RL models.

We find that, after the first $10^4$ training epochs, the RL models optimized using PPO outperform those that were optimized using A2C (Fig. \ref{fig: Lunar Lander}B). Over this interval of training time, we find that the trajectories generated by RL models trained using PPO have smaller maximum eigenvalue norm (corresponding to greater stability -- Fig. \ref{fig: Lunar Lander}C) and larger normalized rank of the controllability matrix (Fig. \ref{fig: Lunar Lander}D). We find that, while in the first 1000 epochs the models trained with A2C exhibit an increase in stability and controllability, further training fails to lead to additional increases (Fig. \ref{fig: Lunar Lander}C, D). This is consistent with the fact that models trained with A2C exhibit worse performance than the models trained with PPO 
(Fig. \ref{fig: Lunar Lander}B).

\section{DISCUSSION}
\label{section: Discussion}

The explosion of development in reinforcement learning (RL) methods \cite{mnih2016asynchronous, schulman2017proximal}, capable of being applied to large and complex environments, has led to RL models with remarkable behavior \cite{mnih2015human, silver2018general, degrave2022magnetic}. Fully understanding these models, at the behavioral level, has remained elusive, in large part because quantifying dynamical behavior is challenging. Here, inspired by recent uses of Koopman operator theory \cite{mezic2005spectral, budivsic2012applied} to analyze and compare machine learning models \cite{ostrow2023beyond, akrout2023representations, redman2024identifying, huang2025inputdsa}, as well as complex systems more generally \cite{rowley2009spectral, avila2020data}, we demonstrate its ability to extract information about the behavior of RL models that has meaningful connections to the associated rewards, on several benchmark environments. 

We find that the fit Koopman with control models shed light on how RL model behavior evolves with training. In some cases, we find that the extracted dynamical features can act as hidden progress measures \cite{barak2022hidden}, identifying improvements of the RL model with training that are not seen by the coarse-grained reward function (Fig. \ref{fig: Cart Pole Hidden Progress Measure}). Monitoring the metrics associated with stability and controllability may therefore be useful for guiding decisions related to training (e.g., deciding how long to train a model). 

The Koopman with control framework also enables comparisons of behavior from RL models trained with different optimizers. This allows for identifications of where the behavior is similar (for instance, after $200$ training epochs on $\texttt{CartPole}$, RL models trained with PPO and A2C have similar reward, maximum eigenvalue norm, and normalized rank of the controllability matrix: Fig. \ref{fig: Cart Pole}B--D) and where the behavior is different (for instance, after $6 \cdot 10^3$ training epochs on $\texttt{Acrobot}$, RL models trained with PPO and A2C have similar reward, but have very different normalized rank of the controllability matrix: Fig. \ref{fig: Acrobot}B, D). Using methods specifically developed to compare control systems, including in environments with partially observable state-spaces (e.g., InputDSA \cite{huang2025inputdsa}), can enable a more principled comparison of the behavior of different RL models. 

Finally, we note that, by examining the dynamical behavior of RL models trained across multiple tasks, we can gain insight into the biases of different optimizers. For instance, RL models trained with A2C see rapid changes in maxmimum eigenvalue norm, as compared to RL models trained with PPO (Figs. \ref{fig: Cart Pole}C, \ref{fig: Acrobot}C, and \ref{fig: Lunar Lander}C). In the case of \texttt{Acrobot}, this change in stability leads to an initial advantage for A2C optimized RL models. However, in the case of \texttt{CartPole}, where the change in stability was in the incorrect direction (decreased stability, as compared to the desired increase in stability), this leads to a need for greater training to improve the reward. Greater comparison between the wide variety of RL optimizers may lead to more efficient and robust training. 

\textit{Limitations.} \hspace{1mm} In this work, we consider RL models trained in environments with discrete-actions spaces, deterministic and time-invariant state transitions, and fully observable state-spaces. While these are significantly simpler than many of the complex environments that RL has been applied to, \texttt{CartPole}, \texttt{Acrobot}, and \texttt{LunarLander} are standard baselines in the RL community. Demonstrating that the Koopman with control framework can shed light on the behavior of RL models in these environments represents an important first step in the larger goal of enhancing the interpretability of RL model behavior. 

While we focus on RL models trained with only two different optimizers (PPO and A2C), we note that PPO continues to be a popular approach and A2C is a simplified version of A3C \cite{mnih2016asynchronous}, an optimizer that enabled state-of-the-art performance on Atari games. Thus, PPO and A2C represent a core baseline of RL optimization methods. In addition, that we find PPO consistently develops -- in all three environments tested -- greater stability (or greater transience, in the \texttt{Acrobot} environment) and greater controllability, supports its observed robustness and wide adoption. 

Lastly, we note that we compute the associated Koopman with control models at a small number of training iterations. Performing our analysis over a wider range of training may identify specific changes in stability and controllability that are training epoch dependent. 

\textit{Future directions.} \hspace{1mm} Our results demonstrate that the Koopman with control framework can provide insight into the behavior of RL models trained on ``physical'' control problems. As RL models can be applied to a very broad range of abstract tasks (e.g., playing chess), a natural next question to investigate is whether and how the Koopman with control framework can be applied to understanding the behavior of RL models in such environments. One potential path for doing this is utilizing a richer selection of observable functions, as opposed to using only time-delays \cite{brunton2017chaos, arbabi2017ergodic}, as was done in this work.  

Finally, we note that recent work has shown that Koopman with control approaches can be used to model the activations of recurrent neural networks, including those underlying RL models \cite{ostrow2023beyond, huang2025inputdsa}. Coupling the activation level investigation pioneered by that work and the behavioral level investigation explored in this work could enable a multi-scale understanding of RL. 

\section*{ACKNOWLEDGMENTS}

We thank Jordan Garrett for assistance with developing the code for running and analyzing the RL models. We thank Ann Huang for assistance in formalizing Eq. 
\ref{eq: RL model update}. We thank Igor Mezić, Yannis Kevrekidis, Mitchell Ostrow, Ann Huang, Leo Kozachkov, Jared Markowitz, and Jordan Garrett for useful discussion on Koopman operator theory for studying RL models. This material is based upon work supported by the Air Force Office of Scientific Research under award number FA9550-22-1-0531.

\small
\bibliography{main.bib}
\bibliographystyle{unsrt}

\end{document}